\documentclass[10pt,reqno]{amsart}
\usepackage{latexsym}
\usepackage{graphicx}
\usepackage{fancyhdr,amssymb}

\usepackage{hyperref} 

\title{A Game of Primes}
\author{Raghavendra N Bhat}
\date{September 2022}
\usepackage{graphicx}
\usepackage{listings}
\usepackage{amsmath}
\usepackage{tabularx}
\usepackage[table]{xcolor}
\usepackage[thinlines]{easytable}
\begin{document}
\maketitle
\begin{abstract}
The basis for most of the ideas mentioned in this paper is the theory of cellular automata. A cellular automata contains a regular grid of cells, with each cell having a pre-defined set of finite states.
The initial state is determined at time/state zero. At this point all the cells are assigned their respective starting states. The automata is defined by a set of simple rules that decide the subsequent states of the cells. We aim to create a cellular automata of prime numbers and come up with some axioms, theorems and conjectures for the same.
\end{abstract}
\section{Introduction}

The main topic of this project is a ``Game of Primes'' which is a cellular automata designed to study prime gaps, prime density and other numerical patterns involving primes.
It is inspired by John Conway's Game of Life \cite{gardner} which was originally defined by Conway in 1970 with the following rules:
\begin{enumerate}
\item Each cell has two states : Alive or Dead.
\item Any live cell with fewer than two live neighbors dies, as if by under population.
\item Any live cell with two or three live neighbors lives on to the next generation.
\item Any live cell with more than three live neighbors dies, as if by overpopulation.
\item Any dead cell with exactly three live neighbors becomes a live cell, as if by reproduction.
\end{enumerate}Let us now define a cellular automata, by the name of ``Game of Primes'' (GOPM). Rather than having an automata on a grid of abstract infinite space and dots, let us attempt to bring these principles of life and automation onto a grid of natural numbers. We define GOPM with the following rules :
\begin{enumerate}
\item We start with an $n$ by $n$ grid of natural numbers which are in arithmetic progression (most of the GOPMs that we discuss in this paper will be grids of consecutive natural numbers).
\item The numbers are filled in a snake-like manner, with the smallest number as the top left corner; and the largest number being the bottom right/bottom left corner, based on parity (most of the GOPMs that we discuss in this paper will have 1 as the starting number, i.e. at the top left).
\item Each cell has two states : Excited or Dormant.
\item If a dormant number is surrounded by  a total of three or more primes and/or excited numbers, it gets excited (reproduction).
\item If an excited number is surrounded by a total of four or more primes and/or excited numbers, it goes back to being dormant (death due to over crowding).
\item If an excited number is surrounded by no other excited numbers and/or primes, it goes back to becoming dormant (death due to loneliness).
\item A ``day" or state is a time interval in which the states of all elements are changed based on the above mentioned rules. For example, day 0 is defined to be the start of the game, where all elements are dormant. Day 1 is the subsequent state which could potentially involve excited cells (depending on the conditions that allow the above mentioned rules to thrive).
\item An ``immortal element" is one that stays excited once it gets excited (i.e, it never goes dormant once it is excited).
\item A ``boundary element" is one that is present along the boundary of the grid. Thus, the four corners, the elements along the top and bottom rows, and the elements along the leftmost and rightmost column make up the set of boundary elements.
\item An NGOPM is a GOPM set up that has 1 as the top-left element and contains consecutive natural numbers from 1 to $n^2$, where $n$ is the dimension. Most of the grids we will be looking at are NGOPMs, as mentioned in rule (1) and rule (2).
\item Thus, a typical GOPM set up is defined by the following variables : $G$, the grid with $k$ or $n$ as the dimension. Two other things determine the $G$; the start number and spacing between the numbers (in NGOPMs we have both the starting number and the gap as 1).\\\\
\begin{center}
\begin{TAB}(e,1cm,1cm){|c:c:c:c:c|}{|c:c:c:c:c|}
1 & 2 & 3 & 4 & 5 \\
10 & 9 & 8  & 7 & 6 \\
11  & 12  & 13  & 14  & 15 \\
20  & 19  & 18  & 17  & 16  \\
21  & 22  & 23  & 24  & 25    
\end{TAB}\\
Game of Primes Interface\\
\end{center}
\end{enumerate}

\section{A GOPM Example}
Let us take a look at a 3x3 NGOPM to understand the above mentioned rules. Initially all the numbers are dormant. As an exercise we can make a note of each number and its prime neighbor count.\\\\
For example : The number 1 has three neighbors out of which 2 and 5 are primes; the number 2 has five neighbors out of which 3 and 5 are primes; the number 3 has three neighbors out of which 2 and 5 are primes.\\\\
Thus, we can conclude that the only three numbers that can get excited are 6, 5, 4 since they have three or more prime neighbors (6 has 2, 5 and 7; 5 has 2, 3, and 7; 4 has 2, 3 and 5).
\begin{center}
\begin{tabular}{|c|c|c|}
  \hline
  1 & 2 & 3\\
  \hline
  6 & 5 & 4\\
  \hline
  7 & 8 & 9\\
  \hline
\end{tabular}
\\
On the zeroth day of the NGOPM 3x3 set up, all elements are dormant\\
\end{center}
\begin{center}
\begin{tabular}{|c|c|c|}
  \hline
  1 & 2 & 3\\
  \hline
  \cellcolor{blue!25}6 & \cellcolor{blue!25}5 & \cellcolor{blue!25}4\\
  \hline
  7 & 8 & 9\\
  \hline
\end{tabular}
\\
On the first day of the NGOPM 3x3 set up, all elements that had three or more prime neighbors get excited (6, 5, and 4)\\
\end{center}
\begin{center}
\begin{tabular}{|c|c|c|}
  \hline
  \cellcolor{blue!25}1 & \cellcolor{blue!25}2 & \cellcolor{blue!25}3\\
  \hline
  \cellcolor{blue!25}6 & 5 & \cellcolor{blue!25}4\\
  \hline
  7 & \cellcolor{blue!25}8 & 9\\
  \hline
\end{tabular}
\\
On the second day of the NGOPM 3x3 set up, 5 goes back to becoming dormant because it has more than three prime+excited neighbors. 6 and 4 continue to live another day. 1, 2, 3 get excited owing to having 4 and 6 as excited neighbors\\
\end{center}
\begin{center}
\begin{tabular}{|c|c|c|}
  \hline
  \cellcolor{blue!25}1 & 2 & \cellcolor{blue!25}3\\
  \hline
  6 & \cellcolor{blue!25}5 & 4\\
  \hline
  \cellcolor{blue!25}7 & 8 & \cellcolor{blue!25}9\\
  \hline
\end{tabular}
\\
On the third day of the NGOPM 3x3 set up, 6, 2, and 4 become dormant and 5 gets back to the excited state. 1, 3, 7, and 9 are the corner elements and since they all have at least one prime neighbor, they stay excited\\
\end{center}

\begin{center}
\begin{tabular}{|c|c|c|}
  \hline
  \cellcolor{blue!25}1 & \cellcolor{blue!25}2 & \cellcolor{blue!25}3\\
  \hline
  \cellcolor{blue!25}6 & 5 & \cellcolor{blue!25}4\\
  \hline
  \cellcolor{blue!25}7 & \cellcolor{blue!25}8 & \cellcolor{blue!25}9\\
  \hline
\end{tabular}
\\
The fourth day of the NGOPM 3x3 set up results in 5 going back to becoming dormant owing to over crowding. 2, 4, 8 and 6 get back to being excited because of prime/excited neighbors\\
\end{center}

\begin{center}
\begin{tabular}{|c|c|c|}
  \hline
  \cellcolor{blue!25}1 & 2 & \cellcolor{blue!25}3\\
  \hline
  6 & \cellcolor{blue!25}5 & 4\\
  \hline
  \cellcolor{blue!25}7 & 8 & \cellcolor{blue!25}9\\
  \hline
\end{tabular}
\\
The fifth day of the NGOPM 3x3 set up is same as the third day with the corners continuing their excited states and 5 getting back to its excited state. We have now reached a cycle since days 3(or 5) and 4 repeat infinitely many times\\
\end{center}
\section{Axioms}
\noindent\textbf{Axiom 3.1 :} \textbf{Every game will eventually have a cycle of states.}\\\\
Since we have a finite number of cells to start with, it is only a matter of time before patterns repeat. With or without the rules, since each cell has two possible states (excited or dormant), the total number of unique states possible is $2^n$ for a grid with $n$ cells. Since the rules of the game are generically defined for all days, a day being reached twice implies that the sequence of days following that repeated day would also repeat; giving rise to a cycle.\\
\\
\textbf{Axiom 3.2 :} \textbf{There has to exist at least one number in the grid that has at least three prime neighbors for the game to start.}\\\\
Rule 4 of the Game of Primes states the requirement of at least three excited/prime neighbors for a dormant cell to get excited. At the start of the game, no cells are excited. Thus, the enablers for excitement are primes only. Hence the axiom.\\
\section{Theorems}
\noindent\textbf{Theorem 4.1 : In an already run GOPM, consider $S$ to be the set of all unique days that have occurred. The maximum number of ``days" that can lead to a particular day in $S$ is two.}\begin{proof} Let $S$ be the set of all unique days in a GOPM setup. Let us define a function from S to S. For every day $x$ $\in$ $S$, $f(x) = y$ where $y$ is the next day. Notice that $y$ $\in$ $S$ as well.\\\\
Assume for the sake of contradiction $\exists$ $a, b, c$ $\in$ $S$ not all equal such that $f(a)$ = $f(b)$ = $f(c)$ = k in a particular GOPM set up $G$. Thus, k follows $a$ as well as $b$. Assume that the sequence of days are $x_1$, $x_2$, $x_3$, $\dots$,$x_n$, $a$, $k$, $y_1$, $y_2$,$\dots$, $y_m, b, k...$By Axiom 1, this leads to a cycle starting at $b$. Thus, the set of all days for G are $\{x_1, x_2, x_3,\dots, x_n, a, k, y_1, y_2, \dots, y_m, b\}$. Thus, $c=a$ or $c=b$. We are done.\\\\
It is important to understand why we are talking about an \textbf{already run} GOPM. If the GOPM has not yet reached a cycle, it is quite possible that there are multiple days that could all lead up to $k$. However, in reality only two of those days end up being a part of the final set $S$ (because the moment the second day leads to $k$, the chain of events repeat). This ensures that the other ``alternate realities" do not occur.\\\end{proof}
\noindent\textbf{Theorem 4.2 : If a GOPM grid has at least one number which has four or more prime neighbors, the cycle length of the GOPM is even.}\begin{proof} Let the cell which has more than three prime neighbors be $x$. By Rule 4 of GOPM, $x$ gets excited as soon as the game starts, i.e. on the first day. It goes back to becoming dormant on the second day owing to Rule 5. By Rule 4, it gets excited once again on the third day. Thus, it alternates between excited and dormant for the rest of the game.\\\\
Assume cycle has been reached. Thus, $\exists$ days $A$ and $B$ such that f(A) = f(B) = $K$ where $K$ is another state (we use the same function defined in Theorem 1). Thus, the state of $x$ in $A$ equals the state of $x$ in $B$ and is different from the state of $x$ in $K$. Hence, the number of days in the cycle $\{K, y_1, y_2...., y_m, B\}$ is even owing to the alternate fluctuation of the state of $x$.\\
\end{proof}
\noindent\textbf{Connections to the Twin Prime conjecture and Prime Gaps:} Suppose there exists a cell $x$ in an NGOPM grid $G$ of dimension $k$ such that $x+1$, $x-1$, $x+k$ and $x-k$ are primes, the cycle length of $G$ is even (the four above mentioned cells are neighbors of $x$). Notice that for this to be possible, the dimension $k$ has to be an odd number and $x$ has to lie in the middle column of the grid (owing to the way the numbers are filled in a GOP grid). Given below is a code in Python which takes the dimension of the grid as an input and prints the smallest possible $x$ (if it exists) that has neighbors as mentioned above. The existence of such an $x$ is a sufficient condition for the cycle length of a grid to be even (not necessary).
\begin{lstlisting}
#Code Starts
def primetest(x):
    if x<0 or x==0 or x==1:
        return False
    y=2
    z=0
    while y**2<=x:
        if x%y==0:
            z=1
            return False
        else:
            y+=1
    if z!=1:
        return True

k=int(input("Enter dimension of NGOPM grid "))
x=k+2
while x<=k**2-k-1:
    if primetest(x+1) == True \
    and primetest(x-1) == True \
    and primetest(x+k) == True \
    and primetest(x-k) == True:
        print(str(x))
        break
    else:
        x+=1
        if x>k**2-k-1:
            print('Did not find x')
#Code Ends
\end{lstlisting}
While this might not be the best way to check for the existence of a cell with at least four prime neighbors, running this code on all odd dimensions from 3-100 (since 1x1 and 2x2 do not have elements which have four neighbors), tells us that 32 of those have such a number $x$.
\begin{center}
\begin{tabular}{|c|c|}
  \hline
  Grid Dimension & Element $x$\\
  \hline
  5 X 5 & 12 \\
  \hline
  7 X 7 & 12\\
  \hline
  11 X 11 & 18\\
  \hline
  13 X 13 & 18\\
  \hline
  17 X 17 & 30\\
  \hline
  19 X 19 & 42\\
  \hline
  23 X 23 & 30\\
  \hline
  29 X 29 & 42\\
  \hline
  31 X 31 & 42\\
  \hline
  37 X 37 & 42\\
  \hline
  41 X 41 & 60\\
  \hline
  43 X 43 & 60\\
  \hline
  47 X 47 & 60\\
  \hline
  53 X 53 & 60\\
  \hline
  59 X 59 & 72\\
  \hline
  61 X 61 & 102\\
  \hline
  67 X 67 & 72\\
  \hline
  71 X 71 & 102\\
  \hline
  73 X 73 & 240\\
  \hline
  79 X 79 & 102\\
  \hline
  83 X 89 & 150\\
  \hline
  89 X 89 & 102\\
  \hline
  97 X 97 & 102\\
  \hline
  101 X 101 & 138\\
  \hline
\end{tabular}
\\The first few prime dimensions which contain an element $x$ with four prime neighbors. We conjecture that this will work for all primes greater than 3.
\end{center}
Interestingly, except for 2 and 3, every prime dimension seems to contain such an element $x$. Iterating this code over primes tells us that all prime dimensions $p$ upto a million have even cycle lengths (owing to the existence of our required cell $x$). Thus, we make the following claim \textbf{``Every prime number $p$, has a number $x$, such that $p+1<=x<=p^2-p-1$ such that $x+1, x-1, x+p$ and $x-p$ are all primes.}\\
\\
If this is indeed true, \textbf{there exist infinitely many NGOPM grids that have even length cycles}.
\\
\\
\textbf{Theorem 4.3 :} \textbf{If a corner element has at least one prime neighbor, it will stay immortal once it gets excited.}\begin{proof} For an excited cell to go back to becoming dormant, there could be either over-crowding or loneliness. Since a corner element has exactly three neighbors, over-crowding is ruled out. Thus, if at least one neighbor is prime, the corner cell remains excited since its time of excitement, as loneliness is ruled out too.\\
\end{proof}
\noindent\textbf{Theorem 4.4 :} \textbf{The following are some Elementary Consequences of the Prime Number Theorem (PNT) :}\\
\begin{enumerate}
\item In an even dimensional NGOPM, the probability that a randomly chosen, non-boundary odd number $o$ gets excited on day 1 is zero, when $o$ does not have 2 as one of its eight neighbors (this does not dependent on PNT and is a result of parity).
\item In an odd dimensional NGOPM, the probability that a randomly chosen, non-boundary even number $n$ gets excited on day 1 is the same as the probability of a randomly chosen non-boundary odd number $m$.
\end{enumerate}
\begin{proof}Theorem 4.4 (2) is a direct consequence of the density of primes. For any given natural number $n$ the number of primes smaller than $n$ is asymptotic to $n/$log $ n$ as per the prime number theorem. Although prime numbers are not uniformly distributed, since our goal is to make a probabilistic argument, we can assume that the probability of randomly picking a prime number is $\frac{\pi(n)}{n}$. Here $\pi$ (n) is the prime counting function. Thus, we can approximate the probability for a number $k<n$ to be prime as 1/log $ n$.
\\\\The main goal is to show that in an odd dimension grid, the probability for any arbitrary element to get excited on day 1 remains approximately the same. Each non-boundary element in an odd dimensional grid has 4 odd neighbors. Since we need at least 3 prime neighbors, the probability of getting excited is asymptotic to $4p^3(1-p) + p^4$, where $p$ is the probability for the odd number to be prime (i.e. $2n/log$ $ n$). However, in an even dimensional grid, an even number has a higher chance of getting excited than an odd number (except for boundary elements and those that have 2 as a neighbor), owing to the fact that it can have up to 6 odd neighbors and an odd number has only 2. Thus, the possibility of a randomly chosen, non-boundary odd number that does not have 2 as a neighbor getting excited is zero on the first day of an even dimensional NGOPM.\\
\end{proof}
\section{Conjectures}
\noindent\textbf{Conjecture 5.1 : For any non-trivial pattern of states of excited elements, there exist infinitely many grids that reproduce the exact pattern of excited states.}
\\\\
This is a direct analogy of the Prime $k$-tuples conjecture \cite{tuth} which states that ``every admissible pattern for a prime constellation occurs infinitely often". The entire structure of the Game of Primes depends upon the starting day, the top left number, and the gap. Thus, if the distribution of primes on the first day can be reproduced (with different numbers but having the exact same number of primes including their distribution and gaps), the entire game will be reproduced. Thus, our goal will be to find two numbers such that, if put at the top of our GOPM grid, will produce the same game. The Prime $k$-tuples conjecture implies there are infinitely many sets of numbers that can duplicate non-trivial prime distributions.\\\\
For example, if we look at numbers from 51-75 (in a 5x5 GOP grid), the primes are 53, 59, 61, 67, 71, 73. Thus, the gaps are 3, 6, 2, 6, 4, 2 (assuming 50 to be equivalent to 0 and 51 being the first number in the grid i.e., top left element and 75 to be the last). This implies that there would be six prime numbers in a GOPM created with dimension 5 and starting number 51.\\\\
The exact same sequence of gaps and number of primes can be achieved in the interval 261-285 (the primes being 263, 269, 271, 277, 281, 283). Hence, starting the grid from 51 or 261 in a 5x5 set up would result in the exact same GOPM. 1281, 14541 and 75981 are some other numbers that can be the starting points in a 5x5 GOPM grid resulting in the same game as 51. Going by the Prime k-tuples conjectures, we should have infinite such numbers that have 6 prime numbers ahead of them with gaps 3, 6, 2, 6, 4, 2.\\\\
Here are the first four ``days" of the 51 and 261 5x5 grids.
We see them to be identical. This shows us that the nature of the primes and their distribution determines the GOPM and not the primes themselves.\\\\
\begin{center}
\begin{tabular}{|c|c|c|c|c|}
  \hline
  51 & 52 & 53 & 54 & 55\\
  \hline
  60 & 59 & 58 & 57 & 56\\
  \hline
  61 & 62 & 63 & 64 & 65\\
  \hline
  70 & 69 & 68 & 67 & 66\\
  \hline
  71 & 72 & 73 & 74 & 75\\
  \hline
\end{tabular}
Day 0
\end{center}
\begin{center}
\begin{tabular}{|c|c|c|c|c|}
  \hline
  51 & 52 & 53 & 54 & 55\\
  \hline
  60 & 59 & 58 & 57 & 56\\
  \hline
  61 & 62 & 63 & 64 & 65\\
  \hline
  70 & \cellcolor{blue!25}69 & 68 & 67 & 66\\
  \hline
  71 & 72 & 73 & 74 & 75\\
  \hline
\end{tabular}
Day 1
\end{center}
\begin{center}
\begin{tabular}{|c|c|c|c|c|}
  \hline
  51 & 52 & 53 & 54 & 55\\
  \hline
  60 & 59 & 58 & 57 & 56\\
  \hline
  61 & \cellcolor{blue!25}62 & \cellcolor{blue!25}63 & 64 & 65\\
  \hline
  \cellcolor{blue!25}70 & \cellcolor{blue!25}69 & \cellcolor{blue!25}68 & 67 & 66\\
  \hline
  71 & \cellcolor{blue!25}72 & 73 & 74 & 75\\
  \hline
\end{tabular}
Day 2
\end{center}
\begin{center}
\begin{tabular}{|c|c|c|c|c|}
  \hline
  51 & 52 & 53 & 54 & 55\\
  \hline
  \cellcolor{blue!25}60 & \cellcolor{blue!25}59 & \cellcolor{blue!25}58 & 57 & 56\\
  \hline
  \cellcolor{blue!25}61 & 62 & 63 & \cellcolor{blue!25}64 & 65\\
  \hline
  70 & 69 & 68 & \cellcolor{blue!25}67 & 66\\
  \hline
  \cellcolor{blue!25}71 & 72 & \cellcolor{blue!25}73 & \cellcolor{blue!25}74 & 75\\
  \hline
\end{tabular}
Day 3
\\The first four days of 5x5 GOPM grid starting from 51.
\end{center}
\begin{center}
\begin{tabular}{|c|c|c|c|c|}
  \hline
  261 & 262 & 263 & 264 & 265\\
  \hline
  270 & 269 & 268 & 267 & 266\\
  \hline
  271 & 272 & 273 & 274 & 275\\
  \hline
  280 & 279 & 278 & 277 & 276\\
  \hline
  281 & 282 & 283 & 284 & 285\\
  \hline
\end{tabular}
Day 0
\end{center}
\begin{center}
\begin{tabular}{|c|c|c|c|c|}
  \hline
  261 & 262 & 263 & 264 & 265\\
  \hline
  270 & 269 & 268 & 267 & 266\\
  \hline
  271 & 272 & 273 & 274 & 275\\
  \hline
  280 & \cellcolor{blue!25}279 & 278 & 277 & 276\\
  \hline
  281 & 282 & 283 & 284 & 285\\
  \hline
\end{tabular}
Day 1
\end{center}
\begin{center}
\begin{tabular}{|c|c|c|c|c|}
  \hline
  261 & 262 & 263 & 264 & 265\\
  \hline
  270 & 269 & 268 & 267 & 266\\
  \hline
  271 & \cellcolor{blue!25}272 & \cellcolor{blue!25}273 & 274 & 275\\
  \hline
  \cellcolor{blue!25}280 & \cellcolor{blue!25}279 & \cellcolor{blue!25}278 & 277 & 276\\
  \hline
  281 & \cellcolor{blue!25}282 & 283 & 284 & 285\\
  \hline
\end{tabular}
Day 2
\end{center}
\begin{center}
\begin{tabular}{|c|c|c|c|c|}
  \hline
  261 & 262 & 263 & 264 & 265\\
  \hline
  \cellcolor{blue!25}270 & 2\cellcolor{blue!25}69 & \cellcolor{blue!25}268 & 267 & 266\\
  \hline
  \cellcolor{blue!25}271 & 272 & 273 & \cellcolor{blue!25}274 & 275\\
  \hline
  280 & 279 & 278 & \cellcolor{blue!25}277 & 276\\
  \hline
  \cellcolor{blue!25}281 & 282 & \cellcolor{blue!25}283 & \cellcolor{blue!25}284 & 285\\
  \hline
\end{tabular}
Day 3
\\The first four days of 5x5 GOPM grid starting from 261.
\end{center}
\textbf{Conjecture 5.2 : The length of the cycle for grids of consecutive natural numbers is seemingly random.}
\begin{center}
\begin{tabular}{|c|c|}
  \hline
  Dimension of GOP & Cycle Length\\
  \hline
  1 X 1 & 0 \\
  \hline
  2 X 2 & 0\\
  \hline
  3 X 3 & 2\\
  \hline
  4 X 4 & 4\\
  \hline
  5 X 5 & 12\\
  \hline
  6 X 6 & 44\\
  \hline
  7 X 7 & 8\\
  \hline
  8 X 8 & 16\\
  \hline
  9 X 9 & 120\\
  \hline
  10 X 10 & 8\\
  \hline
  12 X 12 & 21384\\
  \hline
  17 X 17 & 360\\
  \hline
  19 X 19 & 24\\
  \hline
\end{tabular}
\\
The cycle lengths of some of the first $N$x$N$ grids filled with natural numbers starting with 1\\
\end{center}
The lengths of the cycles for the first $N$x$N$ grids of natural numbers are listed above. We see that the cycle lengths are neither increasing nor decreasing. This does counter the natural intuition which would assume that the cycle lengths should probably increase as the grid size increases (owing to more number of primes and possible patterns).
Here are a few possible explanations for the non-monotone nature of the cycle lengths:
\begin{enumerate}
\item The number of primes in the grid and/or their distribution decides the time taken for the cycle to start but not the length of the cycle itself. Thus, if day $K$ is the start of the cycle, the number of days in between two consecutive appearances of $K$ is independent of the grid size but the days taken to first reach $K$ depends on the grid structure. In this case, we do see some convincing numerical evidence.\\ The number of days taken for the first appearance of the respective $K$ of each of the first 10 Consecutive Natural Number grids is as follows :
\item The Game is partially randomized as it depends on the gaps and distribution of primes. While there are very good estimates for the number of primes smaller than a given number $N$ (the prime number theorem, for example, which states this to be asymptotic to $N/log $ $N$ \cite{zagier}\cite{davenport}), there is not much known about where those primes exactly occur. Moreover, while the ``snake-like" manner of filling the grid does try to keep neighbors on the number line as neighbors on the grid, it does not do it to the elements that are present at the left end of the rows. The top and bottom neighbors of a cell $k$ are $k-n$ and $k+n$ respectively where $n$ is the dimension of the grid.
\end{enumerate}

\section{A Musical Game of Primes}
Rhythm and music have always had connections with mathematics. The Fibonacci sequence \cite{guruprasad} \cite{bhat}, Magic squares \cite{bhat2} and many other mathematical structures, equations and ideas can be represented and appreciated musically. The Game of Primes gives us an interface to create yet another one.\\\\
Count the number of excited cells in each state. Divide it by seven and look at the remainder. For example : If there are 17 cells that are excited, we denote it by 3, 28 excited cells is denoted by 0.\\\\
Thus, we will now have a sequence of numbers ranging from 0-6 for each GOPM with each number denoting the number of excited cells on a particular day modulo 7. Let us now assign the seven fundamental frequency to these numbers from 0-7.\\\\The frequencies that we will consider are 240 Hz, 270 Hz, 288 Hz, 320 Hz, 360 Hz, 405 Hz, 432 Hz for 0-6 respectively. As a side note, these are the frequency notes corresponding to ``Do Re Mi Fa So La Ti".\\\\
For example, looking at the set of days in the 3x3 GOP set up (Section 2), Day 3 = Day 5. Thus, in that case we have a cycle of length 2 starting on the 3rd day. If we were to assign musical notes as discussed above, we would get the following sequence of numbers for days 1-5 with each number being matched to the respective fundamental frequency note :\\
\\
Day 1 : 3 : Fa\\
Day 2 : 6 : Ti\\
Day 3 : 5 : La\\
Day 4 : 1 : Re\\
Day 5 : 5 : La\\
\\
We then have the sequence Re-La repeating indefinitely until the game is stopped.

\section{Future work}
The rules and parameters of the game of primes can be modified to suit our problem of interest. Consider the problem of identifying how ‘fake news’ spreads. By marking the locations of the cells that can create the fake news, one can study what happens to the rest of the cells over a time interval. With appropriate data, it would be possible to add smoothing techniques and reinforcement learning to get better outcomes. The same idea can be used for transportation in countries like India and China that have large populations. With a set of independent variables like number people, roads and peak hours of travel, we hope to tune the dependent variables like number of cars per road, lanes and ride share via cellular automatons.

\section{Compliance with Ethical Standards}
Conflict of Interest: The author declares that there is no conflict of interest.

\section*{Acknowledgment}

The author would like to thank the department of Mathematics at the University of Illinois for their constant support and encouragement. The author would like to mention Professors Alexandru Zaharescu, Bruce Reznick and Bruce Reznick for their inputs. The author extends a special thanks to his father, Dr. Narasimha Bhat for creating an online Game of Primes interface that was of tremendous help for exploring the conjectures and theorems of this paper.

\end{document}